\documentclass[11.5pt]{article}
\usepackage[english]{babel}
\usepackage{amssymb}
\usepackage{graphicx}
\usepackage{graphics}
\usepackage{amsthm}

\input pictex.tex

\newcommand{\clo}{\mathrm{S}^1}

\theoremstyle{definition}

\newtheorem{thm}{Theorem}[section]

\newtheorem{lem}[thm]{Lemma}

\input epsf

\begin{document}
\author{Andr\'es Navas \& Michele Triestino}
\date{}

\title{On the invariant distributions of $C^2$ circle diffeomorphisms 
of irrational rotation number}
\maketitle

Although invariant measures are a fundamental tool in Dynamical Systems, very little is known 
about distributions (i.e. linear functionals defined on some space of smooth functions on 
the underlying space) that remain invariant under a dynamics. Perhaps the most general 
definite result in this direction is the remarkable theorem of A. Avila and A. Kocsard \cite{AK} 
according to which no $C^{\infty}$ circle diffeomorphism of irrational rotation number has 
an invariant distribution different from (a scalar multiple of 
integration with respect to) the (unique) invariant 
(probability) measure. The main result of this Note is an analogous result in low regularity. 
Unlike \cite{AK} which involves very hard computations, our approach is more conceptual. 
It relies on the work of R. Douady and J.-C. Yoccoz concerning automorphic 
measures for circle diffeomorphisms \cite{yoccoz}.

\vspace{0.3cm} 

\noindent{\bf Theorem A.} {\em Circle diffeomorphisms of irrational rotation number 
that belong to the Denjoy class $C^{1+bv}$ have no invariant 1-distributions 
different from the invariant measure.}

\vspace{0.3cm}

Here and in what follows, for $k \geq 0$, by a $k$-distribution we mean a (continuous) linear 
functional defined on the space of $C^k$ real-valued functions on the circle. Notice that for all 
$k' > k$, a $k$-distribution may be seen as a $k'$-distribution, but the converse is false in 
general.

Theorem A allows to show an improved version of the Denjoy-Koksma inequality for $C^1$ test 
functions, thus extending Corollary C of \cite{AK} (valid for diffeomorphisms of class $C^{11}$) 
to $C^{1+bv}$ diffeomorphisms. We omit the proof since it follows the very same lines of that 
of \cite{AK}. Indeed, the only new tool needed in \cite{AK} was the absence of invariant 
1-distributions other than the invariant measure.

\vspace{0.3cm}

\noindent{\bf Corollary.} {\em Let $f$ be a $C^{1+bv}$ circle diffeomorphism of irrational rotation number 
$\rho$. If $(p_k/q_k)$ is the sequence of rational approximations of $\rho$, then for every $C^1$ function 
$u$ on the circle, we have the uniform convergence}
$$u + u \circ f + u \circ f^2 +  \ldots + u \circ f^{q_k - 1} - q_k \int_{\clo} u \hspace{0.04cm} d\mu \longrightarrow 0,$$
{\em where $\mu$ denotes the (unique) invariant probability measure of $f$.}

\vspace{0.3cm}

Letting $u := \log(Df)$ whenever $f$ is a $C^2$ diffeomorphism, this yields a well-known 
result of M. Herman: the sequence $(f^{q_k})$ converges to the identity in 
the $C^1$ topology (see \cite[Chapitre VII]{herman}).

Finally, Theorem A is sharp in that diffeomorphisms with lower regularity may admit invariant 
1-distributions.

\vspace{0.3cm}

\noindent{\bf Theorem B.} {\em For each irrational angle $\rho$, there exists a $C^1$ 
circle diffeomorphism of  rotation number $\rho$ having an invariant 1-distribution 
different from the invariant measure. Moreover, such a diffeomorphism can be taken 
either being minimal or admitting a minimal invariant Cantor set.}

\vspace{0.4cm}

\noindent{\bf I. No invariant distributions for $C^{1+bv}$ minimal diffeomorphisms.}  
Given a $C^{1+bv}$ circle diffeomorphism $f$ of irrational rotation number, let 
$\mu$ be the invariant probability measure of $f$. In order to prove that $f$ has no 
invariant 1-distribution other than $\mu$, it suffices to show that for every $C^1$ 
function $u$ on the circle of zero $\mu$-mean, there exists a sequence of $C^1$ 
functions $v_n$ such that the coboundaries
$$v_n \circ f - v_n$$
converge to $u$ in the $C^1$ topology.\footnote{Actually, a standard application of the 
Hahn-Banach theorem shows that this condition is also necessary, but we will not need 
this fact.} Indeed, if such a sequence exists, then for every invariant 1-distribution $L$,
$$L(u) = \lim_{n \to \infty} \big[ L(v_n \circ f) - L(v_n) \big] = 0.$$

Let us denote by $\lambda$ the (normalized) Lebesgue measure on the circle. The 
existence of the desired sequence $v_n$ is an almost direct consequence of 
the next Proposition.

\vspace{0.42cm}

\noindent{\bf Proposition.}  {\em There exists a sequence of continuous functions $w_n$ such that} 
$$(w_n \circ f ) Df - w_n$$
{\em uniformly converges to $u'$ and for all $n$,}
\begin{equation}\label{zero}
\int_{\clo} w_n d \lambda = 0.
\end{equation} 
\label{key}

\vspace{0.2cm}

Assume this Proposition for a while. 
By integration from $0$ to $x$, we obtain that the function $u(x) - u(0)$ is 
approximated by the sequence of continuous functions 
$$x \mapsto \int_0^x (w_n \circ f ) Df - \int_0^x w_n = 
\int_0^{f(x)} w_n - \int_0^x w_n - \int_0^{f(0)} w_n.$$
We let 
$$v_n(x) := \int_0^x w_n d\lambda,$$
which is well defined due to (\ref{zero}). Then the function 
$u$ is $C^1$ approximated by the sequence $v_n \circ f - v_n + c_n$, 
where $c_n$ is a constant:
$$c_n := u(0) - \int_0^{f(0)} w_n.$$
By integration with respect to $\mu$, one concludes that $c_n$ necessarily 
converges to $0$, thus yielding the desired approximating sequence and 
hence proving Theorem A.

\vspace{0.2cm}

For the proof of the Proposition above, the next particular case (corresponding to the 
displacement function $u(x) := f(x)-x$ minus its rotation number $\rho(f)$) will be crucial.

\vspace{0.42cm}

\noindent{\bf Lemma.}  {\em There exists a sequence 
of continuous functions $\hat{w}_k$ such that} 
$$(\hat{w}_k \circ f ) Df - \hat{w}_k$$
{\em uniformly converges to $Df(x)-1$ and for all $k$,}
$$\int_{\clo} \hat{w}_k d \lambda = 0.$$

\noindent{\bf Proof.} Let 
$$\hat{w}_k := - \frac{1}{q_k} \left[ 1 + Df + Df^2  + \ldots + Df^{q_k-1} \right] + 1,$$
where $(q_k)$ is the sequence of denominators in the rational approximation of 
the rotation number of $f$. Notice that for all $k \geq 1$,
$$\int_{\clo} \hat{w}_k \hspace{0.02cm} d\lambda = 0.$$
Moreover, 
\begin{eqnarray*}
\hat{w}_k (f(x)) Df(x) 
&=& 
- \frac{Df(x)}{q_k} \left[ 1 + Df (f(x)) + Df^2 (f(x)) + \ldots + Df^{q_k-1} (f(x)) \right] + Df(x)\\
&=&
- \frac{1}{q_k} \left[ Df (x) + Df^2 (x) + \ldots + Df^{q_k} (x) \right] + Df(x)\\
&=&
\hat{w}_k (x) - 1 + \frac{1}{q_k} \big[ 1 - Df^{q_k}(x) \big] + Df(x),\\
\end{eqnarray*}
hence
$$Df(x) - 1 = (\hat{w}_k \circ f ) Df - \hat{w}_k + \frac{1}{q_k} \big[ Df^{q_k} - 1 \big].$$
The desired convergence then follows from the Denjoy inequality: 
$$\big| Df^{q_k} \big| \leq \exp(V),$$
where $V$ denotes the total variation of the logarithm of 
$Df$ (see \cite[Chapter 3]{book}). $\hfill\square$

\vspace{0.28cm}

Let us now come back to the Proposition. For the proof, let us recall that 
given $s \in \mathbb{R}$, an $s$-automorphic measure for $f$ is a probability 
measure $\nu$ on the circle such that for every continuous function $\varphi$,
$$\int_{\clo} \varphi \hspace{0.03cm} d\nu 
= \int_{\clo} (\varphi \circ f)  (Df)^s \hspace{0.01cm} d\nu.$$
For a $C^{1+bv}$ circle diffeomorphisms of irrational rotation number, such a measure exists 
and is unique for each $s$ (see \cite{yoccoz}).  The (unique) 
1-automorphic measure is hence the Lebesgue measure. 
Moreover, the uniqueness holds (up to a scalar factor) even 
in the context of signed finite measures (i.e. linear functionals 
defined on the space of continuous functions).

Now, since 
$$\int_{\clo} u' \hspace{0.02cm} d\lambda = 0,$$
a standard application of the Hahn-Banach theorem (to the functional $w \mapsto 
(w \circ f) Df - w$) yields a sequence of continuous functions $\bar{w}_n$ such that 
$$(\bar{w}_n \circ f) Df - \bar{w}_n$$
uniformly converges to $u'$. Indeed, assume otherwise. Then $u'$ does not belong to 
the closure of the set of functions of the form $(w \circ f) Df - w$. The latter set being 
convex, there exists a linear functional $L$ which is identically zero on this set but 
$L(u') = 1$. The former condition means that, as a signed measure, $L$ is 1-automorphic. 
It is hence a nonzero multiple of the Lebesgue measure,  which is absurd since 
$L(u') = 1$ and $u'$ has zero mean with respect to $\lambda$.

The problem with the approximating sequence $\bar{w}_n$ above is that it is unclear 
whether one can ensure that 
$$c_n : = \int_{\clo} \bar{w}_n  d\lambda$$
equals zero. To solve this problem, we consider the functions 
$\tilde{w}_n : = \bar{w}_n - c_n$, which obviously have zero mean (with respect 
to $\lambda$). Given $\varepsilon > 0$, we may choose $n$ such that the absolute 
value of the difference between 
$$(\tilde{w}_n \circ f) Df - \tilde{w}_n$$
and 
$$u' - c_n (Df - 1)$$
is smaller than or equal to $\varepsilon/2$. Moreover, 
the Lemma yields $k = k_n$ such that 
$$c_n \big| (Df - 1) - [(\hat{w}_k \circ f) Df - \hat{w}_k] \big| \leq \frac{\varepsilon}{2}.$$
Putting all of this together, we get that 
$$\big| u' - [(\tilde {w}_n \circ f - c_n \hat{w}_k \circ f ) Df - 
(\tilde {w}_n - c_n \hat{w}_k)\big| \leq \varepsilon,$$
which together with 
$$\int_{\clo} (\tilde {w}_n - c_n \hat{w}_k) \hspace{0.025cm} d\lambda = 0$$
proves the Proposition. 

\vspace{0.4cm}

\noindent{\bf II. Examples of $C^1$ diffeomorphisms with invariant distributions.} 
We first deal with diffeomorphisms with a minimal invariant Cantor set. To do this, 
we next recall a (particular case of a) construction from \cite[Section 5]{yoccoz} and 
then show how this provides an example of a $C^{1}$ circle diffeomorphism of irrational 
rotation number admitting invariant 1-distributions different from the invariant measure.

\vspace{0.35cm}

\noindent{\bf Theorem [R. Douady \& J.-C. Yoccoz].} {\em For each irrational angle $\rho$,  
there exists a $C^1$ circle diffeomorphism $f$ which is a Denjoy counter-example and 
satisfies the following property: for a certain point $x_0$ belonging to the complement 
of the exceptional minimal Cantor set, one has}
$$S := \sum_{n \in \mathbb{Z}} Df^n (x_0) < \infty.$$
{\em In particular, the (probability) measure}
$$\nu := \frac{1}{S} \sum_{n \in \mathbb{Z}} Df^n(x_0) \hspace{0.06cm} \delta_{f^n(x_0)} $$
\noindent{\em is 1-automorphic.}

\vspace{0.35cm}

Let us fix $f$ as above, and let $\mu$ be its (unique) invariant measure. Let us consider the 
linear functional (compare \cite{AK2})
$$L \!: u \mapsto \int_{\clo} u' \hspace{0.01cm} d \nu,$$
defined on the space of $C^1$ functions $u$ on the circle. 
Then $L$ is $f$-invariant. Indeed, 
$$L(u\circ f) 
= \int_{\clo} (u \circ f)' \hspace{0.02cm} d\nu 
= \int_{\clo} (u' \circ f) Df \hspace{0.03cm} d\nu 
= \int_{\clo} u' \hspace{0.02cm} d\nu = L(u),$$
where the third equality follows from that $\nu$ is 1-automorphic. 

To see that $L$ is different from (a multiple of the integration with respect to the) 
invariant measure $\mu$, we let $I$ be the connected component of the complement 
of the exceptional minimal set $K$ that contains $x_0$. Then $K$ coincides with 
the support of $\mu$, so that for every function $u$ supported on $I$, we have 
$$\int_{\clo} u \hspace{0.03cm} d\mu = 0.$$
For such a function, one has 
$$L(u) = \int_{\clo} u' \hspace{0.02cm} d\nu = \frac{u' (x_0)}{S}.$$
However, this expression can take any real 
value for different functions $u$ as above. 

\vspace{0.2cm}

The construction of $C^1$ minimal diffeomorphisms with invariant distributions follows 
a similar strategy. The main tool is a recent result from \cite{matsumoto} which in a 
certain sense can be considered as a measurable counterpart of that of Douady and 
Yoccoz. 

\vspace{0.34cm}

\noindent{\bf Theorem [H. Kodama \& S. Matsumoto].} {\em For each irrational rotation angle 
$\rho$ there exists a minimal $C^1$ circle diffeomorphism $f$ of rotation number $\rho$ 
admitting a measurable fundamental domain $C$ (that is, a measurable set that is 
disjoint from all of its iterates). Moreover, for each point $x_0 \in C$, one has} 
$$S := \sum_{n \in \mathbb{Z}} Df^n (x_0) < \infty.$$

\vspace{0.34cm}

As before, the last condition implies that the measure
$$\nu := \frac{1}{S} \sum_{n \in \mathbb{Z}} Df^n(x_0) \hspace{0.06cm} \delta_{f^n(x_0)}$$
is 1-automorphic for $f$, so that we may again define the $f$-invariant 1-distribution 
$$L\!: u \mapsto \int_{\clo} u' \hspace{0.02cm} d\nu$$
for $C^1$ functions $u$ on the circle. 
However, unlike the previous case, showing that $L$ is nontrivial on the kernel of $\mu$ 
is not geometrically obvious. Nevertheless, there is a simple general argument that shows 
this fact (which also applies to the previous case), as we next explain.

By simple integration, every continuous function on the circle of zero $\lambda$-mean 
can be seen as the derivative of a $C^1$ function having zero $\mu$-mean. 
As a consequence, if $L$ were 
trivial on the kernel of $\mu$, then for every continuous function $v$ on the circle, we would have
$$0 = \int_{\clo} \left( v - \int_{\clo} v \hspace{0.02cm} d\lambda\right) \hspace{0.03cm} d\nu 
= \int_{\clo} v \hspace{0.03cm} d\nu - \int_{\clo} v \hspace{0.03cm} d\lambda.$$
This would mean that $\nu = \lambda$, which is absurd.

\vspace{0.4cm}

\noindent{\bf Remark.} Although Douady-Yoccoz' examples can not be made $C^{1 + \tau}$ 
for any $\tau > 0$, the Denjoy counter-examples extensively described in \cite{book} carry 
1-automorphic atomic measures in the same way as those of \cite{yoccoz}. The case of 
minimal diffeomorphisms is more complicated. Indeed,  Kodama-Matsumoto's examples 
are neither $C^{1+\tau}$, and it is unclear whether one can produce 1-automorphic 
measures by following a similar construction. (Actually, there is no known 
minimal $C^{1+\tau}$ circle diffeomorphism which is non ergodic with 
respect to the Lebesgue measure.) This raises the question of the 
existence of minimal $C^{1+\tau}$ circle diffeomophisms  of irrational 
rotation number having invariant 1-distributions other than the invariant measure.

\vspace{0.4cm} 

\noindent{\bf Question.} Does there exist a $C^2$ circle diffeomorphism of irrational rotation 
number carrying a 2-invariant distribution different from the invariant measure ?

\vspace{0.5cm}

\noindent{\bf Acknowledgments.} The first named author is indebted to A. Kocsard 
for his interest on this Note as well as many useful conversations on the subject. He 
would also like to acknowledge the support of the Fondecyt Grant 1120131 and the  
``Center of Dynamical Systems and Related Fields'' (DySyRF). Both authors would 
like to thank ICTP-Trieste for the hospitality at the origin of this work.


\begin{small}

\vspace{0.2cm}

\noindent Andr\'es Navas\\

\noindent Univ. de Santiago de Chile\\

\noindent Alameda 3363, Estaci\'on Central, Santiago, Chile\\

\noindent E-mail: andres.navas@usach.cl\\

\vspace{0.4cm}

\noindent Michele Triestino\\

\noindent \'Ecole Normale Sup\'erieure de Lyon, CNRS UMR 5669\\

\noindent 46 all\'ee d'Italie, 69364 Lyon 07, France\\

\noindent E-mail: michele.triestino@ens-lyon.fr\\

\end{small}


\begin{thebibliography}{Dillo 83}

\bibitem{AK} {\sc A. Avila \& A. Kocsard.} Cohomological equations and invariant distributions 
for minimal circle diffeomorphisms. {\em Duke Math. J.} {\bf 158} (2011), 501-536.

\bibitem{AK2} {\sc A. Avila \& A. Kocsard.} Invariant distributions for higher dimensional 
quasiperiodic maps. Preprint (2012). 
 

\bibitem{yoccoz} {\sc R. Douady \& J.-C. Yoccoz.} Nombre de rotation des diff\'eomorphismes 
du cercle et mesures automorphes. {\em Regular and Chaotic Dynamics} {\bf 4} (1999), 2-24.

\bibitem{herman} {\sc M. Herman.} Sur la conjugaison diff\'erentiable des diff\'eomorphismes
du cercle \`a des rotations. {\em Publ. Math. de l'IH\'ES}  {\bf 49}  (1979), 5-233.

\bibitem{matsumoto} {\sc H. Kodama \& S. Matsumoto.} Minimal $C^1$-diffeomorphisms 
of the circle which admit measurable fundamental domains. To appear in {\em Proc. of the AMS}.

\bibitem{book} {\sc A. Navas.} {\em Groups of circle diffeomorphisms.}
Chicago Lectures in Mathematics (2011).

\end{thebibliography}
\end{document}